# A Genetic Algorithm Approach for Solving a Flexible Job Shop Scheduling Problem

Sayedmohammadreza Vaghefinezhad[1], Kuan Yew Wong[2]

[1] Department of Manufacturing & Industrial Engineering, Faculty of Mechanical Engineering, Universiti Teknologi Malaysia, 81310 UTM Johor Bahru, Malaysia

[2] Department of Manufacturing & Industrial Engineering, Faculty of Mechanical Engineering, Universiti Teknologi Malaysia, 81310 UTM Johor Bahru, Malaysia

**Abstract**

Flexible job shop scheduling has been noticed as an effective manufacturing system to cope with rapid development in today's competitive environment. Flexible job shop scheduling problem (FJSSP) is known as a NP-hard problem in the field of optimization. Considering the dynamic state of the real world makes this problem more and more complicated. Most studies in the field of FJSSP have only focused on minimizing the total makespan. In this paper, a mathematical model for FJSSP has been developed. The objective function is maximizing the total profit while meeting some constraints. Time-varying raw material costs and selling prices and dissimilar demands for each period, have been considered to decrease gaps between reality and the model. A manufacturer that produces various parts of gas valves has been used as a case study. Its scheduling problem for multi-part, multi-period, and multi-operation with parallel machines has been solved by using genetic algorithm (GA). The best obtained answer determines the economic amount of production by different machines that belong to predefined operations for each part to satisfy customer demand in each period.

**Keywords:** *Flexible Job-Shop Scheduling, Optimization, Flexible Manufacturing System, Integer Programming, Genetic Algorithm.*

## 1. Introduction

Scheduling is significantly investigated in manufacturing systems. Scheduling has been applied to meet customer demand, which plays an important role in customer satisfaction. Additionally, providing a better schedule influences a system's performance. "Sequencing" and "Scheduling" are usually considered as a synonym. Conway, Maxwell, and Miller [1] stated that sequencing problem occurs when there is an order in the number of tasks that must be performed. Baker [2] defined sequencing problem as a particular scheduling problem in which a series of tasks entirely establishes a schedule. The sequence is defined as a permutation of N tasks on specified machines [3]. The job-shop scheduling problem (JSSP) is a general scheduling problem in the manufacturing industry. A classical JSSP consists of N and M number of different jobs and machines, respectively. A set of operations and related machines with known processing time are involved in the process of each job. JSSP is a sequencing problem with no machine substitute for each operation while in FJSSP, alternative machines are considered for each operation. Due to the NP-hard class of FJSSP, meta-heuristic approaches have been widely utilized to solve it. Pezzella, Morganti, and Ciaschetti [4] proposed a genetic algorithm (GA) for solving FJSSP and proved that GA can solve the problem more effectively than tabu search. A GA-based heuristic algorithm has been utilized to solve an integrated scheduling problem consisting of job shop, flow shop and production line [5]. Özgüven, Özbakır, and Yavuz [6] defined routing flexibility as the existence of alternative machines for each operation and process plan flexibility as the existence substitute process plans for each job and developed a mixed-integer linear programming model for FJSSP with routing and process plan flexibility. Fattahi and Fallahi [7] developed a multi objective mathematical model for dynamic FJSSP and solved it using a GA-based algorithm. Li and Chen [8] suggested a new chromosome structure with two rows for machine distribution and working sequence and designed the relevant GA operators. Zhang, et al. [9] improved GA to generate a better initial population and suggested an enhanced chromosome structure for solving FJSSP more conveniently. Savaş [10] suggested a GA with novel operations in order to find the minimum makespan in a scheduling problem with dissimilar parallel machines. Chen, Wu, Chen, and Chen, [11] proposed FJSSP with parallel machines and used GA to select machine and sequence of operations for each machine. Majority of previous studies tried to minimize the makespan by selecting the appropriate machine for each operation. In this paper, the objective function is maximizing the total profit by:

- ✓ Determining the amount of production in normal time and overtime.

The remainder of this paper are structured as follows: Section 2: Description of proposed model, Section 3: Genetic algorithm, Section 4: Case study, Section 5: Result's interpretation, Section 6: Concluding remarks.

## 2. Problem Description

In this paper, an integer mathematical model is proposed to solve FJSSP over a finite time horizon. The objective function is to maximize total profit. Several operations must be performed in order to complete a particular part. Parallel machines have been defined for each operation. All predefined machines for each operation can participate in production. Furthermore, by assuming time-varying costs of raw materials and fluctuating selling prices of completed parts, the optimum answer satisfies the customer demand, by holding the economic amount of inventory. Moreover, overproduction will be sold less than the selling price at the end of the time horizon.

### 2.1 Assumptions of the model

To solve this problem, some assumptions have been made:

- Parts have numerous operations, which have to be completed before shipment to the final store.
- Alternative machines have been considered for each operation.
- Customer demand, selling price and cost of each raw material are not constant and can be estimated in the beginning of time horizon.
- Machines' normal time and overtime are constant and out sourcing is not done.
- There is no back order.
- The holding cost is constant at all time.
- Demand can be fulfilled with inventory from the previous period.
- General purpose machines are assumed.
- There is no waste for each operation. It means that all non-conforming items can be reworked.
- Work-in-process inventory is not allowed. It means that parts can only be held when all operations in their processes have been finished.

### 2.2 Indices

The below indices have been defined for model formulation.

i     index for part types ($i = 1,\ldots,P$)

j     index for machine types ($j = 1,\ldots,M$)

k     index for number of operation for each part ($k = 1,\ldots,K_i$)

t     index for periods ($t = 1,\ldots,T$)

### 2.3 Input Parameters

The required input parameters are:

$D_{it}$    Demand for part i in period t
$C_{ikj}$    Cost of the $k^{th}$ operation for part i by using the machine type j in normal time ($/min)
$C^+_{ikj}$    Cost of the $k^{th}$ operation of part i by using the machine type j in overtime ($/min)
$S_{it}$    Selling price for part i in period t ($)
$S'_i$    Selling price for part i at the end of time horizon ($)
$P_{ikj}$    Process time of the $k^{th}$ operation of part i by using machine type j (min)
$B_{jt}$    Normal time capacity of machine j in period t (min/month)
$B^+_{jt}$    Overtime capacity of machine j in period t (min/month)
$W_i$    Weight of part i (kg)
$RP_{it}$    Raw material cost for part I in period t ($/kg)
$H_i$    Holding cost for part i ($/unit per period)

$$Y_{ikj} = \begin{cases} 1 & \text{if machine j is assigned to the } k^{th} \text{ operation for part i} \\ 0 & \text{otherwise,} \end{cases}$$

### 2.4 Decision variables

The decision variables of this mathematical model have been defined as follow:

$X_{ikjt}$    Amount of part i that will be created after the $k^{th}$ operation of its process route by machine j in period t and in normal time.

$X'_{ikjt}$    Amount of part i that will be created after the $k^{th}$ operation of its process route by machine j in period t in overtime.

### 2.5 Mathematical model

By considering the above notations, the created model is:

$$Max\ z = \sum_{i=1}^{P}\sum_{t=1}^{T} D_{it}\ S_{it}$$
$$+ \sum_{i=1}^{P}\sum_{t=1}^{T}\sum_{k=1}^{K_i}\sum_{j=1}^{M} ([(X_{ikjt} + X'_{ikjt})/K_i] - D_{it})\ S'_i$$
$$- \sum_{i=1}^{P}\sum_{k=1}^{K_i}\sum_{j=1}^{M}\sum_{t=1}^{T} X_{ikjt}\ Y_{ikj}\ P_{ikj}\ C_{ikj}$$
$$- \sum_{i=1}^{P}\sum_{k=1}^{K_i}\sum_{j=1}^{M}\sum_{t=1}^{T} X'_{ikjt}\ Y_{ikj}\ P_{ikj}\ C^{+}_{ikj}$$
$$- \sum_{i=1}^{P}\sum_{k=1}^{K_i}\sum_{j=1}^{M}\sum_{t=1}^{T} [(X_{ikjt} + X'_{ikjt})/K_i]\ W_i\ RP_{it}$$
$$- \sum_{i=1}^{P}\sum_{t=1}^{T}\sum_{k=1}^{K_i}\sum_{j=1}^{M} ([(X_{ikjt} + X'_{ikjt})/K_i] - D_{it})\ H_i \quad (1)$$

S.t.

$$\sum_{j=1}^{M} X_{ikjt} + X'_{ikjt} = \sum_{j=1}^{M} X_{i(k+1)jt} + X'_{i(k+1)jt}$$
$$i \in \{1, 2... P\},\ k \in \{1, 2... K_i-1\},\ t \in \{1, 2... T\} \quad (2)$$

$$\sum_{j=1}^{M}\sum_{t'=1}^{t}(X_{ikjt'} + X'_{ikjt'})\ Y_{ikj} \geq D_{it}$$
$$i \in \{1, 2... P\},\ k \in \{1, 2... K_i\},\ t \in \{1, 2... T\} \quad (3)$$

$$\sum_{i=1}^{P}\sum_{k=1}^{Ki} X_{ikjt}\ P_{ikj}\ Y_{ikj} \leq B_{jt}$$
$$j \in \{1, 2... M\},\ t \in \{1, 2... T\} \quad (4)$$

$$\sum_{i=1}^{P}\sum_{k=1}^{Ki} X'_{ikjt}\ P_{ikj}\ Y_{ikj} \leq B^{+}_{jt}$$
$$j \in \{1, 2... M\},\ t \in \{1, 2... T\} \quad (5)$$

$X_{ikjt}$, $X'_{ikjt}$   Integer

$i = 1, 2... P,\ j = 1,2...M,\ k = 1, 2... K_i,\ t = 1, 2...T$ (6)

In the above model, Eq. (1) is the objective function, which represents the total profit. The first and second terms calculate the gross profit of sales in each period and at the end of the time period, respectively. The third and fourth terms compute the total operation cost in normal and overtime. The fifth term determines the total cost of raw materials while the last term in the objective function calculates the holding cost. Since, the amount of each part is transferred between the operations of part's process, in the second, fifth and sixth terms of objective function, the summation of these amounts has been divided by the number of operations for each part. Eq. (2) provides a restriction that the amount of each part remains constant in all operations. Eq. (3) indicates that the amount of parts produced must satisfy the demand. Eq. (4) and Eq. (5) consider the machine capacity limitations during normal time and overtime. Eq. (6) indicates that $X_{ikjt}$ and $X'_{ikjt}$ are integer decision variables.

## 3. Genetic Algorithm

This kind of problem is known as NP-hard problem that can be solved more efficiently with GA and other metaheuristic methods. In this study, GA was selected to solve it. The major mechanisms of GA are: 1. Creating the chromosome representation, 2. Creating the primary population, 3. Creating the adjustment function for fitness evaluation, 4. Defining the selection strategy, 5. Selecting genetic operators for creating a new generation, 6. Defining the parameter values

3.1 Chromosome Structure

An array of direct values, has been used to form a chromosome with the length of sum ($K_i$ * T) for each working state (normal/overtime), since the total length of a chromosome is 2 * sum ($K_i$ * T). There is also a set of sub-chromosomes for each of them that is related to alternative machines with the length of number of available machines for performing each operation. The value of genes indicates the amount of production. Fig. 1 shows an example for the structure of the proposed chromosome for one part and one period when there are two operations for completing this part, and machine M1 is related to the first operation and machines M3 and M6 are assigned to the second operation.

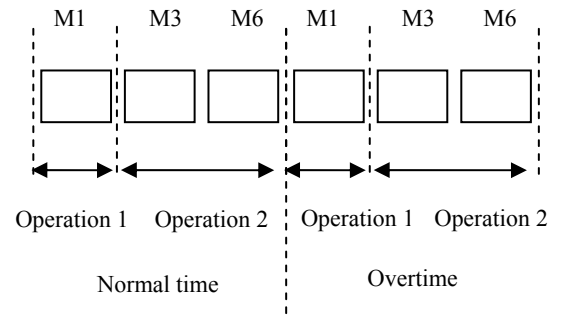

Fig. 1 The Structure of Proposed Chromosome

3.2 Selection Strategy

The tournament method has been utilized as the selection strategy in this paper. This method randomly selects the number of individuals equal to the tournament size, and then chooses the best individual out of that set to be a parent.

## 3.3 Genetic Operators

The following procedure describes the crossover operation of GA: The algorithm generates two random integers m and n between 1 and the length of the chromosome and chooses genes numbered less than or equal to m from the first parent, selects genes numbered from m+1 to n from the second parent, and selects genes numbered greater than n from the first parent. The algorithm then merges these genes to form a single chromosome.

The mutation is performed in two steps. First, the algorithm generates a random number between 0 and 1 and assigns them to genes of each chromosome. Wherever this number is less than or equal to the mutation probability, the related gene will be chosen for mutation. In the second step, the algorithm substitutes each selected gene by a random number selected uniformly from the range.

## 3.4 Termination Criteria

The algorithm will be stopped if it reaches a specified maximum number of generations or if it reaches a specified maximum number of iterations without any improvement.

## 4. Case Study

The selected case study produces some parts of industrial valves. Table 1 shows the parts' weight and demand in the second and third columns. Table 2 and Table 3 show the selling price of part and cost of raw materials that have been forecasted for three periods. Table 4 shows the operations of parts' processes and substitute machines for each operation. The last column demonstrates the related process time for these alternative machines. The machines' available time and associated operation cost are presented in Table 5. Furthermore, the holding cost is 0.1 dollar per month per unit.

Table 1: Weight and Demand for each Part

| Part | $W_{i\,(Kg)}$ | $D_{it}$ | | |
|---|---|---|---|---|
| | | $t = 1$ | $t = 2$ | $t = 3$ |
| Part 1 | 0.168 | 4200 | 4500 | 4300 |
| Part 2 | 0.207 | 3500 | 2500 | 2750 |
| Part 3 | 0.5 | 3000 | 2800 | 3000 |

Table 2: Selling Price for each Part

| Part | $S_{it}$ ($) | | | $S'_i$($) |
|---|---|---|---|---|
| | $t = 1$ | $t = 2$ | $t = 3$ | |
| Part 1 | 1.6 | 1.65 | 1.65 | 0.206 |
| Part 2 | 1.7 | 1.75 | 1.7 | 0.279 |
| Part 3 | 2.98 | 3 | 3.1 | 0.675 |

Table 3: Raw Material Cost for each Part

| Part | $RP_{it}$ ($) | | |
|---|---|---|---|
| | $t = 1$ | $t = 2$ | $t = 3$ |
| Part 1 | 2.23 | 2.35 | 2.45 |
| Part 2 | 2.5 | 2.5 | 2.7 |
| Part 3 | 2.6 | 2.6 | 2.7 |

Table 4: List of Alternative Machines and Related Process Time for the Part's Operation

| Part (i) | Operation(k) | Machine(j) | $P_{ikj}$(min) |
|---|---|---|---|
| Part 1 | 1 | 1,2,3 | 0.5, 0.5, 0.3 |
| Part 1 | 2 | 4,5,6,7 | 1.2, 1.4, 1.5, 1 |
| Part 1 | 3 | 8,9 | 1.5, 1 |
| Part 2 | 1 | 4,5,6,7 | 1.3, 1.5, 1.6, 1.1 |
| Part 2 | 2 | 8,9 | 2.5, 2 |
| Part 3 | 1 | 8,9 | 1, 2 |
| Part 3 | 2 | 4,5,6,7 | 0.6, 0.8, 0.9, 0.4 |
| Part 3 | 3 | 4,5,6,7 | 0.8, 0.9, 1 ,0.7 |

Table 5: Machine's Available Time and Operation Cost in Normal time and Overtime

| Machine | Operation Type | $B_{jt}$ (min/month) | $B^+_{jt}$ (min/month) | $C_{ikj}$ ($/min) | $C^+_{ikj}$ ($/min) |
|---|---|---|---|---|---|
| M1 | Shot Blast | 9240 | 2700 | 0.1 | 0.15 |
| M2 | | 9240 | 2700 | 0.1 | 0.15 |
| M3 | | 9240 | 2700 | 0.12 | 0.18 |
| M4 | CNC | 21600 | 5280 | 0.3 | 0.45 |
| M5 | | 21600 | 5280 | 0.25 | 0.375 |
| M6 | | 21600 | 5280 | 0.2 | 0.3 |
| M7 | | 21600 | 5280 | 0.33 | 0.495 |
| M8 | Assembly | 21600 | 5280 | 0.05 | 0.075 |
| M9 | | 21600 | 5280 | 0.08 | 0.12 |

## 5. Computational Results

The proposed problem has been solved using Matlab software version R2011b and implemented in an Intel Core i5 CPU 2.53 GHz with a 6 GB memory and Windows 7-based operation system. The results based on 0.9 crossover rate, 0.05 mutation rate and 5,000 generations have been collected in Table 6 to Table 8. The objective function was $24862.69. Table 6 shows the obtained answer in the first period. The fourth column indicates the amount of production for each part by different machines of each operation in normal time while the fifth column reveals the amount of production in overtime. The next column shows the summation of production by different machines of each operation for each part in normal time and overtime. Since it has been

defined in Eq. (2), these amounts remain equal in different operations for each part. The next column shows the part demand in the first period while the last column shows the amount of overproduction in the first period. The overproduction amounts must be held for future demands of the next period and they are shown as $I_{it}$ in Table 7 and Table 8. Table 7 and Table 8 demonstrate the final solution for the second and third periods.

Table 6: Final Results for the First Period

| Part (i) | Operation (k) | Machine (j) | $X_{ikjt}$ (t = 1) | $X'_{ikjt}$ (t = 1) | $X_{ikjt} + X'_{ikjt}$ | $D_{it}$ | $(X_{ikjt} + X'_{ikjt}) - D_{it}$ |
|---|---|---|---|---|---|---|---|
| Part 1 | 1 | 1,2,3 | (1304), (1404), (654) | (511), (430), (414) | 4717 | | |
| Part 1 | 2 | 4,5,6,7 | (1587), (1004), (489), (555) | (211), (157), (352), (362) | 4717 | 4200 | 517 |
| Part 1 | 3 | 8,9 | (1377),(1340) | (707),(1293) | 4717 | | |
| Part 2 | 1 | 4,5,6,7 | (765), (461), (434), (695) | (158), (433), (137), (417) | 3500 | 3500 | 0 |
| Part 2 | 2 | 8,9 | (1463),(1531) | (280),(226) | 3500 | | |
| Part 3 | 1 | 8,9 | (1312),(1468) | (269),(111) | 3160 | | |
| Part 3 | 2 | 4,5,6,7 | (767), (442), (266), (230) | (730), (122), (135), (468) | 3160 | 3000 | 160 |
| Part 3 | 3 | 4,5,6,7 | (608), (949),(506), (297) | (417), (210), (151), (22) | 3160 | | |

Table 7: Final Results for the Second Period

| Part (i) | Operation (k) | Machine (j) | $X_{ikjt}$ (t = 2) | $X'_{ikjt}$ (t = 2) | $X_{ikjt} + X'_{ikjt}$ | $I_{i1}$ | $D_{it}$ | $I_{i2}$ |
|---|---|---|---|---|---|---|---|---|
| Part 1 | 1 | 1,2,3 | (613), (601), (931) | (621), (674), (723) | 4163 | | | |
| Part 1 | 2 | 4,5,6,7 | (428), (744), (986), (388) | (892), (405), (4), (316) | 4163 | 517 | 4500 | 180 |
| Part 1 | 3 | 8,9 | (1088),(1355) | (370),(1350) | 4163 | | | |
| Part 2 | 1 | 4,5,6,7 | (650), (50), (85), (201) | (278), (584), (403), (249) | 2500 | 0 | 2500 | 0 |
| Part 2 | 2 | 8,9 | (1560),(766) | (89),(85) | 2500 | | | |
| Part 3 | 1 | 8,9 | (721),(1151) | (790),(27) | 2689 | | | |
| Part 3 | 2 | 4,5,6,7 | (119), (432), (10), (496) | (274), (245), (504), (609) | 2689 | 160 | 2800 | 49 |
| Part 3 | 3 | 4,5,6,7 | (425), (737), (769), (245) | (196), (47), (58), (212) | 2689 | | | |

Table 8: Final Results for the Third Period

| Part (i) | Operation (k) | Machine (j) | $X_{ikjt}$ (t = 3) | $X'_{ikjt}$ (t = 3) | $X_{ikjt} + X'_{ikjt}$ | $I_{i2}$ | $D_{it}$ | $I_{i3}$ |
|---|---|---|---|---|---|---|---|---|
| Part 1 | 1 | 1,2,3 | (146), (891), (357) | (967), (730), (1031) | 4122 | | | |
| Part 1 | 2 | 4,5,6,7 | (1037), (582), (187), (432) | (543), (788), (365), (188) | 4122 | 180 | 4300 | 0 |
| Part 1 | 3 | 8,9 | (1431),(926) | (1003),(762) | 4122 | | | |
| Part 2 | 1 | 4,5,6,7 | (518), (219), (41), (500) | (109), (537), (646), (183) | 2753 | 0 | 2750 | 3 |
| Part 2 | 2 | 8,9 | (1363),(1317) | (73),(0) | 2753 | | | |
| Part 3 | 1 | 8,9 | (936),(1492) | (27),(496) | 2951 | | | |
| Part 3 | 2 | 4,5,6,7 | (165), (351), (1306), (66) | (38), (247), (470), (308) | 2951 | 49 | 3000 | 0 |
| Part 3 | 3 | 4,5,6,7 | (598), (133), (161), (73) | (450), (340), (156), (1040) | 2951 | | | |

## 6. Conclusions

In this paper, an integer mathematical model has been developed to solve FJSSP with time-varying demands, fluctuating raw material costs and different selling prices over the finite time planning horizon. The objective function was maximizing the total profit while meeting some constraints. A GA solver was programmed using Matlab software version R2011b to solve this model. A computational study has been conducted for an industrial case company which produces some parts of industrial gas valves. The obtained results were acceptable and all constraints have been satisfied.

### Acknowledgment

Authors would like to appreciate the Case Study Company and Universiti Teknologi Malaysia for their support and assistance.

## References

[1] RW. Conway, WL. Maxwell, and LW. Miller, Theory of Scheduling, Massachusetts: Addison-Wesley Publishing Company, 1967.

[2] KR. Baker, Introduction to Sequencing and Scheduling, New York: John Wiley & Sons, 1974.

[3] M. Pinedo, Scheduling: Theory, algorithms, and systems, 2nd ed. New Jersey: Prentice-Hall Inc., 2002.

[4] F. Pezzella, G. Morganti, and G. Ciaschetti, "A genetic algorithm for the Flexible Job-shop Scheduling Problem," Computers & Operations Research, vol. 35, no. 10, 2008, pp. 3202-3212.

[5] H. Yan and X. Zhang, "A Case Study on Integrated Production Planning and Scheduling in a Three-Stage Manufacturing System," IEEE Transactions on Automation Science and Engineering, vol. 4, no. 1, 2007, pp. 86-92.

[6] C. Özgüven, L. Özbakır, and Y. Yavuz, "Mathematical models for job-shop scheduling problems with routing and process plan flexibility," Applied Mathematical Modelling, vol. 34, no. 6, 2010, pp. 1539-1548.

[7] P. Fattahi and A. Fallahi, "Dynamic scheduling in flexible job shop systems by considering simultaneously efficiency and stability," CIRP Journal of Manufacturing Science and Technology, vol. 2, no. 2, 2010, pp. 114-123.

[8] Y. Li and Y. Chen, "A Genetic Algorithm for Job-Shop Scheduling," Journal of Software, vol. 5, no. 3, 2010, pp. 269–274.

[9] G. Zhang, L. Gao, Y. Shi, "An effective genetic algorithm for the flexible job-shop scheduling problem," Expert Systems with Applications, vol. 38, 2011, pp. 3563-3573.

[10] B. Savaş, "Non-identical parallel machine scheduling using genetic algorithm," Expert Systems with Applications, vol. 38, 2011, pp. 6814-6821.

[11] J. C. Chen, J.-J. Wu, C.-W. Chen, and K.-H. Chen, "Flexible job shop scheduling with parallel machines using Genetic Algorithm and Grouping Genetic Algorithm," Expert Systems with Applications, 2012.

**Sayedmohammadreza Vaghefinezhad** received Master of Engineering degree (Industrial Engineering) from UTM. He graduated in B.S course in the field of Industrial Engineering (Industrial Production) at Sharif University of Technology. His current research interests include operations research, simulation of operations, production and operations management, modeling and analysis of operation systems, evolutionary algorithms, and computer programming languages.

**Kuan Yew Wong** holds a PhD from the University of Birmingham, England. Currently, he heads the Industrial Engineering Laboratory at Universiti Teknologi Malaysia (UTM). His research interests are related to operations management and industrial engineering.